\documentclass[a4paper, 11pt,notitlepage]{report}

\usepackage[english]{babel}
\usepackage{amsmath, amssymb,amsthm}
\usepackage[utf8]{inputenc}
\usepackage{graphicx}
\usepackage{enumerate}
\usepackage{dsfont}
\usepackage{slashed}
\usepackage{textcomp}
\usepackage{url}
\usepackage{color}
\usepackage{placeins}
\usepackage{multirow}
\usepackage{subfigure}
\usepackage{braket}
\usepackage{mathtools}

\usepackage{titling}

\usepackage[hidelinks]{hyperref}
\usepackage{csquotes}
\usepackage{algorithm}
\usepackage{algorithmic}

\def \foldfix {}
\def \foldpicpap {}

\newcommand{\costsim}{c_{sim}}
\newcommand{\jacob}{J}

\newcommand{\timperiod}{\theta}
\newcommand{\idn}{\mathds{1}}
\newcommand{\nrm}[1]{\lVert {#1} \rVert}
\newcommand{\elrmat}[3]{#3 \in \mathbb{R}^{#1 \times #2}}
\newcommand{\cels}{\,^\circ\text{C}}

\newcommand{\nDAEsolv}{\# solves DAE}

\newcommand{\cons}{G}
\newcommand{\consmn}{G_0 }
\newcommand{\consmax}{\bar G}
\newcommand{\daytmp}{T_{d}}
\newcommand{\feedin}{P}
\newcommand{\pinlim}{\bar \feedin}
\newcommand{\minediff}{\bar \edens}

\newcommand{\constrt}{K}

\newcommand{\timoptcon}{\mathcal{T}_c}
\newcommand{\timopt}{\mathcal{T}_d}

\newcommand{\violcon}{\nrm{u-\tilde u}_2/ \nrm{\tilde u}_2}
\newcommand{\violfeed}{\nrm{P-\pinlim}_{\text{max}}/\pinlim}
\newcommand{\pumppow}{P_{\text{hyd}}}
\newcommand{\Tminhs}{T_{\text{min}}^{\text{cons}}}
\newcommand{\Tmaxnet}{T^{\text{net}}_{\text{max}}}

\newcommand{\dt}{\text{dt}}

\newcommand{\ematph}{Q}
\newcommand{\wgtfun}{\gamma}
\newcommand{\nvflow}{n_F^v}
\newcommand{\nqflow}{n_F^q}

\newcommand{\edginN}{\delta^{\text{in}}(N)}
\newcommand{\edgoutN}{\delta^{\text{out}}(N)}
\newcommand{\pipes}{\mathcal{P}}
\newcommand{\npip}{n_P}
\newcommand{\pathhs}{K}

\newcommand{\disref}{{FOMR}}
\newcommand{\timlev}{\mu}
\newcommand{\maxreldev}{\max_h{\nrm{y^h-y^h_R}_2/\nrm{y^h_R}_2}}

\newcommand{\maxpressdiff}{\overline{\Delta p}}

\newcommand{\maxpreshous}{p_\text{max}^{\text{h}}}
\newcommand{\minpreshous}{p_\text{min}^{\text{h}}}
\newcommand{\presdifhs}{\Delta p^h}
\newcommand{\minovhs}{\min_{h \in \mathcal{H}}}
\newcommand{\maxovhs}{\max_{h \in \mathcal{H}}}

\newcommand{\retdens}{\edens_R}
\newcommand{\edens}{e}
\newcommand{\param}{\kappa}
\newcommand{\uT}{u_T}
\newcommand{\up}{u_p}
\newcommand{\parcont}{u_T^\param}
\newcommand{\indepvol}{\tilde q}

\newcommand{\pipind}{\alpha}
\newcommand{\cellind}{\beta}
\newcommand{\mempip}[1]{{#1}^\pipind}

\newcommand{\Rbody}{\mathds{R}}
\newcommand{\algreq}{\REQUIRE}

\newcommand{\algstat}{\STATE}
\newcommand{\algwhil}{\WHILE}
\newcommand{\algednwhil}{\ENDWHILE}

\newtheorem{proposition}{Proposition}

\usepackage[affil-it]{authblk}

\title{Optimal control of district heating networks using a reduced order model}

\author[1,2]{Markus Rein}
\author[2]{Jan Mohring}
\author[1]{Tobias Damm}
\author[1]{Axel Klar}
\affil[1]{\footnotesize TU Kaiserslautern, 67663 Kaiserslautern, Germany}
\affil[2]{\footnotesize Fraunhofer Institute for Industrial Mathematics ITWM, 67663 Kaiserslautern, Germany}

\begin{document}
\maketitle
\begin{abstract}

We study the optimal control of district heating networks using a reduced order model based on a system theoretic description close to the underlying Euler equations. In the presented scenarios, the central task is to limit the maximal feed-in power occurring as a product of control and state variables. The underlying dynamics of heating networks acting as optimization constraints pose the central computational complexity, prohibiting the determination of an optimal control online. The advection of the injected energy density on the network results in an index-1, quadratic in state differential algebraic equation, challenging to reduce. The suggested reduced model decreases the computation time of the optimization significantly. The effectiveness of the presented approach is demonstrated for an existing, large-scale heating network including changes of flux directions.\\

\end{abstract}

\section{Introduction}
\begin{figure}
\centering
\includegraphics[scale=0.4]{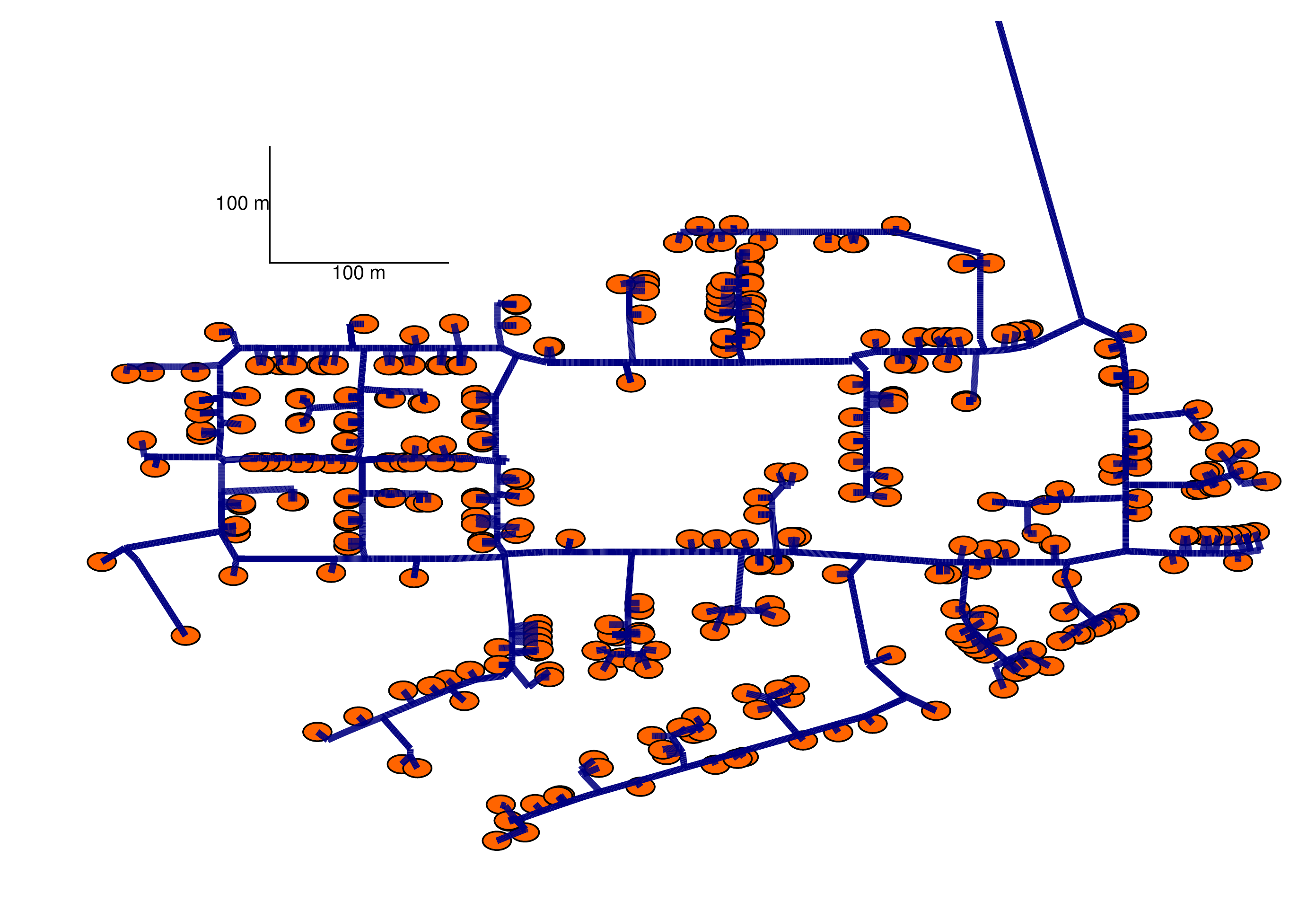}
\caption{Topology of an existing heating network supplying a district. The flow part of the network contains 333 consumers (circles), and 775 pipelines (lines). The data was supplied by Technische Werke Ludwigshafen AG.}
\label{fig_geom_neub}
\end{figure}

In this contribution an optimal control problem for district heating networks is solved using reduced order models. Heating networks are of particular interest for low-carbon energy supply due to their flexibility in using different sources of energy \cite{rezaie_district_2012, talebi_review_2016}. The energy density $\uT$ injected at a power plant is guided to consumers of different sizes using a network of pipelines referred to as flow network. At the consumers, the local volume flow is regulated using heat exchangers to match the time dependent power consumption $\cons$ given the currently available energy density $\edens$. Fig.~\ref{fig_geom_neub} illustrates an existing large scale network considered in this contribution. Its outline data is supplied in tab.~\ref{tab:neub}. A central aim of operating these networks lies in efficiently planning the input energy density $\uT$. It defines the power feed-in $\feedin=(\uT-e_R) \hat q$ in combination with the aggregated volume flow $\hat q$, and the energy density of the cooled fluid $e_R$ entering the plant in the return network. Due to the high transport times from source to consumption in large scale networks, the power feed-in $\feedin(t)$ at each point of time $t$  needs not to match the current power consumption $\cons$ at the heat exchangers. While both quantities are coupled by conservation of energy within the heating network, there exists an essential optimization potential in distributing $\uT$ over time. Injecting a high energy density in times of large volume flows $\hat q$ requires firing additional vessels which might be unfavorable for economic and ecological reasons. Vice versa, planning might also enable to use external overcapacities of energy resulting from renewable energies. In the specific application, a waste to heat incineration plant is able to deliver power up to a maximum level at no cost. When exceeding this level for only a short time interval, significant costs result from using gas boilers to cover peaks in the injected power.\\

In the optimization task, the dynamical energy transport on these complex networks is an essential constraint to ensure a physically relevant control and requires a significant amount of computation time. Since controls of power plants are updated every $15\,$min in application, the corresponding optimization needs to be sufficiently fast. To this end, the usage of reduced order models is central \cite{Benner2014PDE, hovland_explicit_2008, lass_model_2017}. For network systems, both tasks of defining an optimal control \cite{colombo_optimal_2009}, as well as formulating reduced order models are complex. Due to their mathematical complexity and the large benefits, the optimization of energy networks such as electric \cite{gottlich_electric_2016}, water \cite{geisler_mixed_2011}, gas and heating networks is an active field of research. Concerning the model order reduction of energy systems, different works already exist for gas networks \cite{egger_structure-preserving_2018, Prog_DAE_2015}, electric networks \cite{benner_model_2011} and other applications \cite{borggaard2015model}. Similarly, many publications focus on the optimization of gas networks \cite{mehrmann_model_2018,geisler_solving_2015,geisler_mixed_2011,manchanda_challenges_2017,gugat_mip-based_2018,gugat_towards_2018}. While gas networks are modeled by compressible Euler equations, the transport fluid for heating networks is water in the fluid phase inducing incompressible dynamics. On the one hand this simplifies the equations for the conservation of mass and momentum significantly. On the other hand, for heating networks the conservation of inner energy describing the dynamical transport of thermal energy is the dominant effect and mathematically the most challenging one. It leads to a time-dependent advection of the injected thermal energy yielding a large, dynamically changing delay in the energy transport between source and consumption. Towards the formulation of a reduced order model, advection on network systems is a demanding problem.\\ 

Selected works discussing the optimal control of heating networks are mentioned subsequently. In \cite{guelpa_optimal_2016}, pumping costs resulting from the variation of pressure and massflow are optimized using a reduced order model without reflecting the dynamically changing thermal transport. In \cite{sandou_predictive_2005}, a predictive controller is formulated. The advection of thermal energy is reflected by the method of characteristics on each pipeline assuming a constant time delay between source and households. In \cite{baviere_optimal_2018} the transport of thermal energy from source to each consumer is described by virtual single pipelines including a constant time delay from source to consumption points. In this contribution, similar to the approach used for gas networks, a model very close to the underlying Euler equations is formulated, allowing to precisely model the thermal transport dynamics and addressing the central difficulties of heating networks. These are dynamically changing delay times from source to consumers, as well as changes of flux direction. Using a system theoretic approach to model the dynamics allows to use effective tools from model order reduction.
\begin{table}[h!]
\centering
\begin{tabular}{@{}llllll@{}}
\hline
Edges & Nodes & Loops & Pipelines & Consumers  &Pipeline length \\ \hline
1108 & 770 & 6 & 775 & 333 & 8676\,m \\
\hline
\end{tabular}
\caption{Outline data for the flow part of the heating network presented in fig.~\ref{fig_geom_neub}.}
\label{tab:neub}
\end{table}
\\ \\
The article is structured as follows. After introducing the model for the energy transport within heating networks in section \ref{sec:model_dynamics}, a reduced order model required to efficiently simulate the relevant outputs of the heating network, is described in section \ref{sec_ROM}. Hereafter, the optimal control problem is stated in section \ref{sec_cont_prob}, and the numerical determination of an optimal control is discussed in section \ref{sec_det_optcont}. Finally, the application of the reduced order model for determining an optimal control is numerically investigated in section \ref{sec_numeris} and compared to different fidelities of full order models.

\section{Model for the transport dynamics}\label{sec:model_dynamics}
In this section, we suggest a system description of the dynamics of heating networks beneficial for the formulation of a reduced order model. The assumptions used in the derivation are the following. Since water is the transport medium, the incompressible limit is used. Furthermore, the pipelines transporting the fluid are assumed to be perfectly isolated avoiding thermal losses. In agreement with the incompressibility assumption, the density is assumed to be constant. Since for heating networks the pressure difference introduced by friction at the pipeline border dominates the pressure difference induced by acceleration, inertial effects are neglected. Finally, the return temperatures of all heat exchangers are assumed to be equal and constant, in line with the available data. Since these return temperatures define the input for the return network, this allows to neglect the return network when modeling heating networks, because it contains no observable relevant for the dynamics. Under these assumptions, the Euler equations of the transported fluid within a pipeline $\pipind \in \pipes$ read
\begin{eqnarray}
0 &=& v^\pipind_x \label{eq_incomp}\\
0 &=& \mempip{\Delta p} + \frac{\lambda^\pipind \rho \mempip l}{2\mempip d}|v^\pipind |v^\pipind + \rho g \mempip{\Delta z}\label{eq_consmom_sim}\\
\dot \edens^\pipind &=& -v^\pipind \edens^\pipind_x \label{eq_consen_sim},
\end{eqnarray}
where $\pipes$ denotes the set of pipes including $\npip$ elements. While the conservation of mass (\ref{eq_incomp}) degenerates to the incompressibility assumption, the conservation of momentum was integrated along the pipeline length $l^\pipind$ to obtain (\ref{eq_consmom_sim}), introducing pressure and height differences $\Delta p$, $\Delta z$ over pipeline borders. The loss of pressure by friction with the pipeline border is modeled using the Darcy-Weisbach equation \cite{welty_fundamentals_2008}. It couples to a dimensionless friction factor $\lambda^\pipind$ depending on the properties of pipeline $\pipind$. The change of the friction factor with Reynolds number is neglected, allowing to set $\lambda^\pipind$ constant in time within this contribution. The remaining quantities are the pipeline diameter $d^\pipind$, the density $\rho$, and the gravitational constant $g$. Eq.~(\ref{eq_consen_sim}) describes the advection of the energy density $\edens$ subject to the velocity $v$. While the energy density at the start of the pipe is defined by the inflow, the velocity is defined by the consumers stations entering as additional network components. Transferring the pipeline model to the network context, the resulting algebraic constraints at the nodes are given by
\begin{eqnarray}
\sum_{\pipind \in \edgoutN} \mempip q(t) \mempip \edens(t,0) &=& \sum_{\pipind \in \edginN} \mempip q(t) \mempip \edens(t,\mempip l)\label{eq_consen_node}\\
\mempip \edens(t,0) &=& \edens^N(t),\; \pipind \in \edgoutN\label{eq_mix}\\
\sum_{\pipind \in \edgoutN} \mempip q(t) &=& \sum_{\pipind \in \edginN} \mempip q(t)\label{eq_cons_vol}\\
\sum_{\pipind \in O} \frac{\lambda^\pipind \rho l^\pipind}{2d^\pipind}|v^\pipind (t)|v^\pipind (t)  &=& 0, \; O \in \mathcal{L}\label{eq_Kirch2}\\
\big[ q^i(t) \cdot (\edens^i(t,0)-\retdens) \big]_{i \in \mathcal{H}} &=& \cons(t)\label{eq_hous}\\
\edens^s(t,l^s) &=& \uT(t) \label{eq_input_flow}\\
p^s(t,l^s) &=& \up(t) \label{eq_input_press}.
\end{eqnarray}
In the presented network description, $\pipes$ denotes the set of pipelines, $\mathcal{H}$ the set of heat exchangers or consumers. Each element in $\mathcal{L}$ represents a set of paths forming a network loop. The presented algebraic equations describe the conservation of energy (\ref{eq_consen_node}) and volume (\ref{eq_cons_vol}) over node $N$, where $\edginN$ ($\edgoutN$) denotes edges entering (exiting) node $N$. Eq. (\ref{eq_mix}), claims that different incoming energy densities instantly mix to a defined outgoing energy density identical for all edges exiting node $N$. Eq. (\ref{eq_Kirch2}) results from claiming continuity of the pressure, what defines pressure levels at each node $N$. Summing the resulting pressure differences along loop $O$ results in (\ref{eq_Kirch2}). As the second of the two flow defining equations, (\ref{eq_hous}) defines how the power consumption $\cons$ is covered by the difference in energy density at each house and its volume flow. Each of the corresponding heat exchangers ensures by regulation of the volume flow $q$ that each consumer receives the required power consumption $\cons$ for an arbitrary energy difference within the allowed range. Finally, (\ref{eq_input_flow}) and (\ref{eq_input_press}) set the controls for energy density $\uT$ and the pressure level $\up$ at the power plant entering the flow network.\\

For a numerical treatment of the partial differential algebraic equation (PDAE) (\ref{eq_incomp}-\ref{eq_input_flow}), we perform a spatial discretization of the energy density $\edens$ using the upwind scheme \cite{leveque_numerical_2008}. To this end, $\edens_{\pipind, \cellind}$ denotes cell $\cellind$ on pipeline $\pipind \in \pipes$. This allows to formulate the differential algebraic equation (DAE) used to simulate the dynamics of the full order model,
\begin{equation}\label{eq_sys}
\begin{aligned}
\dot \edens &= A(v) \edens + B(v) \uT (t),\\
y &= C \edens , \\
0 &= g(q,y,\cons).
\end{aligned}
\end{equation}
Here $e$ denotes the vector of energy densities of all finite volume cells on all pipelines. The operator $\elrmat{n}{n}{A}$ couples both finite volume cells locally on a pipeline and globally over every node including conservation of energy. The operator $\elrmat{n}{1}{B}$ defines where the thermal input is applied to the network. The output matrix $\elrmat{o}{n}{C}$ measures observables at consumer stations. The system description (\ref{eq_sys}) is Lyapunov stable \cite{rein_model_2019} with known energy matrix $\ematph$. Hence, (\ref{eq_sys}) can be reduced in a stability preserving way as described in section \ref{sec_ROM}. Eq.~(\ref{eq_sys}) forms a DAE of index 1. For the rest of the contribution, we refer to (\ref{eq_sys}) as full order model (FOM). Different full order models resulting from varying number of finite volume cells are evaluated towards their applicability in determining an optimal control in section \ref{sec_numeris}. This includes the comparison to a reduced order model (ROM) formulated in the following section.

\section{Formulation of a reduced order model}\label{sec_ROM}
Subsequently, the formulation of a reduced order model reducing the computational cost of solving the dynamics of heating networks (\ref{eq_sys}) in time is summarized. A detailed derivation of the reduction technique is given in \cite{rein_model_2019}.\\

As shown in \cite{rein_model_2019}, the velocity dependent system operators $A(v), B(v)$ defined in (\ref{eq_sys}) allow for an affine decomposition with respect to the velocity,
\begin{align}\label{eq_affine_decomp}
A(v) &= \sum_{i=1}^{\nvflow} \wgtfun^v_i(v)  A^v_i \\
B(v) &= \sum_{i=1}^{\nvflow} \wgtfun^v_i(v)  B^v_i,
\end{align}
where $B^v_i, A^v_i$ are time-independent matrices multiplied with the velocities defining the time-dependent weighting function $\wgtfun : \Rbody^{\nvflow} \rightarrow \Rbody^{\nvflow}$. Here, $\nvflow$ denotes the number of velocity configurations required to describe the dynamics. The function $\wgtfun_i$ equals the velocity $v^i$ if $A^v_i$ is valid regardless of the flux directions in the network. If in contrast $A^v_i$ is only valid for a specific flow direction, $\wgtfun_i$ is a nonlinear function. If no changes of the flux direction occur in the network, $\wgtfun_i = v^i, \; \forall i \in \pipes$. The determination of all possible flux directions is performed in the offline-phase of the reduction process. This description allows to form a reduced model offline based on projection, where $\elrmat{n}{r}{V,W}$ denote Petrov-Galerkin projection matrices \cite{antoul_approx}. It can be shown that Lyapunov stability translates to the ROM, when implying the energy matrix $\ematph$ in the projection step \cite{polyuga_structure_2010, gugercin_structure-preserving_2012}. For the model of heating networks considered here, the construction of the corresponding energy matrix $\ematph$ is shown \cite{rein_model_2019}. Specifically, volume conservation (\ref{eq_cons_vol}) defined by the current velocity field ensures Lyapunov stability of the FOM \cite{rein_model_2019}. Hence, it is advantageous to rewrite the affine decomposition (\ref{eq_affine_decomp}) in terms of independent volume flows $\indepvol$, what automatically ensures the volume conservation. This is possible, since (\ref{eq_cons_vol}) can be solved a-priori based on the network diameters, forming the set of independent volume flows $v = N \indepvol$. The resulting reduced order model replaces the full order DAE (\ref{eq_opt_DAE}) in finding an optimal control and is given by
\begin{equation}\label{eq_sys_red}
\begin{aligned}
\dot \edens_r &= \sum_{i=1}^{\nqflow} \wgtfun^q_i(\indepvol) W^T A_i^v V   \edens +  \sum_{i=1}^{\nqflow} \wgtfun^q_i(\indepvol) W^T B_i \uT (t)\\
y_r &= C V \edens\\
0 &= g(\indepvol,y_r,\cons).
\end{aligned}
\end{equation}
Since system (\ref{eq_sys_red}) changes dynamically, a global Galerkin projection has to be determined, reflecting all relevant linear submodels $\mathcal{Q}^e$ \cite{benner_survey_2015}. The latter are determined using a greedy-strategy \cite{haasdonk_reduced_2008}. Based on the current projection $V$, the linearization $v^* \in \mathcal{Q}^\delta \supset \mathcal{Q}^e$ is determined, exhibiting the largest error in the transfer function between full and reduced order model. Here, $\mathcal{Q}^\delta$ is the set of linearizations, for which the deviation is tested. It is determined based on training simulations spanning the range of scenarios for which the ROM should be applied hereafter. For the linearization $v^*$, the local Galerkin projection is performed and added to the global projection using a singular value decomposition. This process is repeated, until all linearizations of both full and reduced order model exhibit a relative error of the transfer function below $\Delta$. To obtain the local Galerkin projection, we use a moment-matching technique in frequency space \cite{rein_model_2019}. Similar to IRKA \cite{gugercin_$mathcalh_2$_2008}, it forms a reduced order model interpolating the transfer function of the full order model for a fixed volume flow field at certain points in frequency space. The designed ROM should be valid for different environmental temperatures and only depend on the network topology. Thus, the time in the offline-phase to form the ROM is not taken into account.\\

\section{Control problem}\label{sec_cont_prob}
For the set of discrete points of time $\timopt = \{t_0,...,t_e\}$, and the corresponding continuous interval $\timoptcon = [t_0,t_e]$ we search a parameterized control $\parcont : \timoptcon \rightarrow \Rbody$ of the input energy density, minimizing the following objective 
\begin{align}
J(\parcont) &= \eta_1\nrm{\dot{u}_T^{\param} (\cdot)}_{l_2(\timopt)}^2  +\nrm{u_T^{\param} (\cdot)- \eta_2}_{l_2(\timopt) }^2 \label{eq_objec},
\end{align}
subject to
%
\begin{align}
\parcont(t)&\leq \bar u , \quad t \in \timopt \label{eq_op_uptemp} \\
\quad \edens_{i,n}(t) &\geq \minediff_i , \quad t \in \timopt  \quad \forall i\in \mathcal{H} \label{eq_op_minTemdiff} \\
p^i(t) &\leq \maxpreshous, \quad t \in \timopt \quad i \in \mathcal{H}  \label{eq_op_maxPres_hs}\\
p^i(t) &\geq \minpreshous, \quad t \in \timopt \quad i \in \mathcal{H} \label{eq_op_minPres_hs}\\
(\parcont(t)-\retdens) \sum_{i \in \mathcal{H}} q^i(t) &\leq \pinlim, \quad t \in \timopt \label{eq_op_Pin}  \\
D(\dot \edens(t), \edens(t), \edens_0, q(t), \parcont, \cons(t)) &= 0 \quad t \in \timopt \label{eq_opt_DAE}.
\end{align}
The objective function (\ref{eq_objec}) penalizes the temporal variation and the distance to the regularization parameter $\eta_2$ of the parameterized control $\parcont$. The regularization parameters $\eta_1, \eta_2$ are used to equalize both contributions in the objective which are motivated as follows. Minimizing the temporal variation leads to realistic controls, which can be realized properly by the power plant. By choosing $\eta_2$ sufficiently small, the mean value of the control decreases which systematically reduces thermal losses by cooling effects. Decreasing the injected energy density will lead to higher pressures in the network, which are restricted as explained below. Constraints (\ref{eq_op_uptemp}-\ref{eq_op_minPres_hs}) are technical restrictions in line with standard operation instructions formulated as optimization constraints. The upper energy limit (\ref{eq_op_uptemp}) forces the fluid to be in the liquid phase. Note that energy densities in the network can never exceed the one at the inlet. In (\ref{eq_op_minTemdiff}) a minimal energy density is required for proper operating conditions of the heat exchangers. In addition, the pressure levels of consumption nodes in the flow network are restricted to upper and lower bounds (\ref{eq_op_maxPres_hs},\ref{eq_op_minPres_hs}). Eq.~(\ref{eq_op_Pin}) sets an upper bound to the maximal injected power, which avoids the use of additional energy sources. Below this limit, the power plant can supply demands by energy stemming from a waste to heat incineration plant at no costs. Due to the transport time of the injected energy from source to consumers, the control has a delayed effect on the consumer, allowing to influence the temporal distribution of the injected power. Finally, (\ref{eq_opt_DAE}) reflects the energy transport (\ref{eq_consen_sim}) along the network restricted by the algebraic coupling conditions (\ref{eq_consen_node}-\ref{eq_input_press}) and the initial state $\edens_0(\parcont)$, cf. section \ref{sec:model_dynamics}. For full and reduced order models, (\ref{eq_opt_DAE}) is replaced by (\ref{eq_sys}) and (\ref{eq_sys_red}) respectively.\\

The consumption $\cons$ is assumed to be known a-priori, which is a typical assumption in the simulation of heating networks. Since cooling effects are neglected and the energy densities in the return network are modeled equal and constant, an open loop control problem results. The return network exhibits a constant energy density $\retdens$ entering the feed-in power (\ref{eq_op_Pin}) as a parameter.\\

The required pumping power $\pumppow$ necessary to retain a fluid flow in the network is stemmed by the pumps in the depot. It is bounded above as follows,
\begin{align}\label{eq_pump_pow}
\pumppow(t) = \Delta p^s(t) \sum_{i \in \mathcal{H}} \frac{\cons_i(t) }{\edens_i(t) - \retdens} \leq \frac{\max_{t \in \timoptcon}(\Delta p^s(t)) }{\min_{t \in \timoptcon}{\edens(t)}} \max_{t \in \timoptcon}\sum_{i \in \mathcal{H}} \cons_i(t),
\end{align}
where $\Delta p^s$ describes the pressure difference achieved at the source edge in the depot. For typical networks, the maximum pressure difference at the depot is smaller than $10\,$bar. Approximating the energy density by $\edens \approx \rho c_p T$, with material constants described in section \ref{sec_numeris}, a maximum aggregated power consumption of $1\,$MW leads to a corresponding pumping power of $16\,$kW. Thus the pumping power is suppressed by almost two orders of magnitude compared to the thermal power.
\subsubsection*{Treatment of pressure constraints}
Subsequently, we explain how to fulfill the pressure constraints (\ref{eq_op_maxPres_hs}, \ref{eq_op_minPres_hs}) in a simplified manner. It relies on limiting the difference of the maximum and minimum pressure levels measured at all consumption points. This allows to adjust the pressure control $\up : \timoptcon \rightarrow \Rbody$ after finding the optimal control of the energy density.\\

By (\ref{eq_consmom_sim}), the pressure difference from source to consumption point $h$ is defined by
\begin{align}\label{eq_pressdiffhs}
\presdifhs \equiv p^h - p^s = \rho g (z^s - z^h) - \rho \sum_{i \in \pathhs_h}\lambda^i  \frac{l^i}{2d^i} v^i |v^i| ,
\end{align}
where $\pathhs_h$ denotes an arbitrary path from the source to consumer $h \in \mathcal{H}$. Although the pipeline velocities change dynamically, diameters, lengths and the height profile on the path from source to each household determine the resulting pressure difference to a large extent. This stabilizes the constraint limiting the maximum pressure difference and motivates the following proposition.
%
\begin{proposition}\label{prop_press}
The pressure constraints (\ref{eq_op_maxPres_hs},\ref{eq_op_minPres_hs}) are satisfied by defining an alternative constraint on the difference of the maximum and minimum pressure realized at all consumption points,

\begin{align}\label{eq_limabspress}
\max_{h \in \mathcal{H}}(\tilde p^h(t)) - \min_{h \in \mathcal{H}}(\tilde p^h(t)) \leq \maxpressdiff \leq \maxpreshous -\minpreshous \quad t \in \timopt,
\end{align}
where $\tilde p^h$ denotes the pressure level at consumption point $h$ resulting from a simulation with source pressure $\tilde{u}_p$. Here, $\maxpressdiff$ is the true limit for the pressure difference entering the alternative optimization constraint (\ref{eq_limabspress}). The control of the pressure level leading to admissible pressures at consumption points is obtained a-posteriori by a time dependent shift,
\begin{align}\label{eq_postercont}
\up(t) = \minpreshous - \min_{h \in \mathcal{H}}(\presdifhs(t))\quad  t \in \timopt.
\end{align}
\begin{proof}
By limiting the difference between maximal and minimal absolute pressure levels in (\ref{eq_limabspress}), their value relative to the pressure control is limited as well,
\begin{align}\label{eq_limpressdiff}
\maxovhs(\tilde p^h - \tilde{u}_p) - \min_{h \in \mathcal{H}}(\tilde p^h - \tilde{u}_p) = \max_{h \in \mathcal{H}}(\presdifhs) - \min_{h \in \mathcal{H}}(\presdifhs) \leq \maxpreshous -\minpreshous,
\end{align}
where $\presdifhs$ denotes the pressure difference from house to source, which is independent of the pressure control at the source by (\ref{eq_pressdiffhs}). Thus, the pressure level at each consumption point for an arbitrary source pressure $p^s$ reads,
\begin{align}
p^h = p^s + \presdifhs.
\end{align}
Inserting the suggested control (\ref{eq_postercont}) allows to determine the minimum pressure at each consumption point by
\begin{align}
\minovhs(p^h) = p^s + \minovhs(\presdifhs) = \minpreshous - \minovhs(\presdifhs) + \minovhs(\presdifhs) = \minpreshous.
\end{align}
Similarly, the maximum pressure level is limited by
\begin{equation}\label{eq_estmaxpress}
\begin{aligned}
\maxovhs(p^h) &= p^s + \maxovhs(\presdifhs) \\
&= \minpreshous - \minovhs(\presdifhs) + \maxovhs(\presdifhs)\\
&\leq \minpreshous + \maxpreshous - \minpreshous\\
&=\maxpreshous,
\end{aligned}
\end{equation}
where the inequality in (\ref{eq_estmaxpress}) is obtained by using (\ref{eq_limpressdiff}).
\end{proof}
\end{proposition}
\subsubsection*{Feed-in power and control of energy density}
The feed-in power is the central constraint to limit additional costs, since it avoids the usage of additional energy resources. It is defined by
\begin{align}\label{eq_feedin}
\feedin &= (\parcont(t)-\retdens) \sum_{i \in \mathcal{H}} q^i(t)\\
&= (\parcont(t)-\retdens) \sum_{i \in \mathcal{H}} \frac{\cons_i(t)}{\edens_{i,n_i}(t) - \retdens}.
\end{align}
The control $\parcont$ affects the feed-in $\feedin$ (\ref{eq_feedin}) in two ways. First, by setting the current input energy density $\parcont$ and second, by defining the volume flow which results from the current energy densities at heat exchangers. These in turn equal the control $\parcont(\tau)$ at a past time $\tau$. Depending on the current state $\edens_{h,n_h}$ at consumer stations, the input control can both amplify and weaken the feed-in power with regard to the current consumption $\cons$. In the stationary case $\edens = u_0$, where $u_0$ denotes the constant input, the feed-in power is the temporally shifted consumption profile. Hence, it also matches the high characteristic power peaks in the morning and the evening hours. In contrast, by anticipating the expected consumption and the transport time of the injected power, peaks in the injected power can be reduced. Since the determination of an optimal control is initialized with a constant temperature, the red, solid lines in parts (b) of fig.~\ref{fig_conv_-3}, \ref{fig_conv_3}, \ref{fig_conv_7,5} visualize the consumption profile equaling the displayed feed-in.
\section{Determination of an optimal control}\label{sec_det_optcont}
Subsequently, we discuss the computation of an optimal control for the problem (\ref{eq_objec}-\ref{eq_opt_DAE}). The main idea is to eliminate the transport dynamics from the optimization constraints by solving them explicitly and passing the remaining constraints to the MATLAB nonlinear optimization tool fmincon.\\

Supplying an initial control which satisfies all constraints is an open problem. Hence, the initial parameter set $\param_0$ generally at least violates the feed-in constraint (\ref{eq_op_Pin}) and the determination of a feasible solution is performed initially. For the current parameter vector $\param_i$, the transport dynamics described by the DAE (\ref{eq_opt_DAE}) are solved along $\timopt$. In the solution process, both the trajectory of state variables and their parameter gradients are calculated. This allows to evaluate the true optimization constraints (\ref{eq_op_uptemp}-\ref{eq_op_minTemdiff}, \ref{eq_op_Pin}, \ref{eq_limabspress}), and their gradients with respect to the current parameter $\param_i$. In this step, the limits for the pressure levels (\ref{eq_op_maxPres_hs}, \ref{eq_op_minPres_hs}) are replaced by (\ref{eq_limabspress}) as described in proposition \ref{prop_press}. This allows to focus on the determination of the thermal control $\uT^{\param_i}$ in solving the optimal control problem. The required pumping power resulting from the pressure control $\up$ can be neglected as described in section \ref{sec_cont_prob}. Solving the dynamics explicitly avoids the large computational cost of passing them as optimization constraints. Values and gradients of the optimization constraints for the current parameter are passed to fmincon using the active set method together with the value and parameter of the objective function (\ref{eq_objec}). A summary of the algorithm used to determine the optimal control is provided in alg.~\ref{alg_detoptcont}.\\
\begin{algorithm}
\caption{Numerical computation of an optimal control}
\label{alg_detoptcont}
\begin{algorithmic}[1]
\algreq Initial parameter set $\param_0$, convergence tolerance of nonlinear optimization.
\algwhil{convergence tolerance not satisfied}
\algstat{Solve DAE (\ref{eq_opt_DAE}) using the implicit midpoint rule (\ref{eq_impmid}) for the current parameter vector $\param_i$.
\begin{align*}
D(\dot \edens(t), \edens(t), \edens_0, q(t), \uT^{\param_i}, \cons(t)) =0, \quad  t \in \timopt.
\end{align*}}
\algstat{Determine constraints $\constrt$ defined in (\ref{eq_op_uptemp}-\ref{eq_op_minTemdiff}, \ref{eq_op_Pin}, \ref{eq_limabspress}), and their parameter gradients $\partial_{\param_i}\constrt$ based on the solution of (\ref{eq_opt_DAE}).}
\algstat{Evaluate objective function $J(\uT^{\param_i})$ defined in (\ref{eq_objec}) and its gradient $\partial_{\param_i} J(\uT^{\param_i})$.}
\algstat{Update parameter $\param_{i+1} \leftarrow \text{fmincon}(J(\uT^{\param_i}),\;\partial_{\param_i} J(\uT^{\param_i}),\; \constrt,\; \partial_{\param_i}\constrt)$.}
\algednwhil
\algstat{Adjust pressure control $\up$ according to proposition \ref{prop_press}.}
\end{algorithmic}
\end{algorithm}
%
%
%
\subsubsection*{Extraction of parameter gradients}
To estimate the effect of a change in the parameterized control on the relevant outputs of the heating networks, the sensitivities of both the objective function and the constraints with respect to the parameters have to be determined $\forall t \in \timopt$. To this end, gradients of both the control and state variables regarding the control parameters have to be extracted from the forward solution of the DAE (\ref{eq_opt_DAE}). For the input signals typically applied to heating networks, an implicit time integration of the DAE proved to be beneficial. A general implicit time integration scheme, in which $x_\timlev$ and $x_{\timlev+1}$ denote the state variables at former and future time levels can be written as
\begin{align}
\tilde f(x_\timlev,x_{\timlev+1},u) = 0.
\end{align}
The derivative of the future state variable $\partial_\param x_{\timlev+1}$ is obtained by the derivative of the old state variable $\partial_\param x_{\timlev}$ using the implicit function theorem,
\begin{align}\label{eq_implicit}
\partial_\param x_{\timlev+1} &= -\left(\frac{\partial \tilde f}{\partial x_{\timlev+1}} \right)^{-1} \left( (\partial_{x_{\timlev}} \tilde f) \partial_\param x_\timlev+ \partial_u \tilde f \partial_\param u \right).
\end{align}
Hence, based on the sensitivity of the initial state $\partial_\param \edens_{\pipind, \cellind}(t=0)$, the gradient information can be propagated along the solution of the DAE. Using (\ref{eq_implicit}) allows to determine the gradient \textit{after} solving for the new time step. This is in contrast to many automatic differentiation approaches in which the gradient information has to be tracked during the determination of the future time layer causing additional computational cost.\\

For the time integration of the DAE in this contribution, the implicit midpoint rule as a second order symplectic integrator is used,
\begin{align}\label{eq_impmid}
\tilde f(x_\timlev,x_{\timlev+1},u) = x_{\timlev+1} - x_\timlev - \dt f\left(t_\timlev + \frac{\dt}{2}, \frac{1}{2}(x_\timlev + x_{\timlev+1}),u\right ),
\end{align}
where $\dt$ denotes the time step, and $f$ the dynamical part of the DAE (\ref{eq_opt_DAE}).
\section{Numerical validation}\label{sec_numeris}
\subsection{Time integration of the forward simulation}
The solution of the forward problem (\ref{eq_opt_DAE}) within the determination of an optimal control is performed by the system descriptions (\ref{eq_sys}, \ref{eq_sys_red}). To solve the DAEs (\ref{eq_sys}, \ref{eq_sys_red}) within the time horizon required for the optimal control, the implicit midpoint-rule (\ref{eq_impmid}) is used. Full order models are unreduced ($W=V=\idn$), while for the reduced order model a Galerkin projection is applied. Sparse matrix operations are considered in the full order case. To solve systems (\ref{eq_sys},\ref{eq_sys_red}) efficiently, a domain decomposition is performed \cite{rein_model_2019}. Different parts of the network are treated as independent systems, with their linkages moderated by artificial inputs. This allows to form the system operators and solve the nonlinear equations introduced by the implicit time integration scheme efficiently. Full and reduced order models are simulated using the decomposition. Due to the affine system representation, the determination of the Jacobian can be determined analytically for both full- and reduced order models. The simulations presented in the following sections are performed using MATLAB(R) R2016b on an Intel(R) XEON(R) CPU E5-2670 processor @ 2.60GHz. The nonlinear system of equations resulting from both the algebraic equations and the time integration scheme are solved using the MATLAB function fsolve.

\subsection{Definition of test scenarios}\label{sec_def_scen}
To demonstrate the effectiveness of the reduced order model, the large-scale heating network presented in fig.~\ref{fig_geom_neub} is studied. The loops visible in the left part of the network pose a central difficulty since the thermal transport can take different paths at the same time to reach a certain destination in the network. Moreover, changes of flux direction occur, which change the set of possible paths the transported quantity takes dynamically.\\

The robustness of the reduced model towards its application in the optimization is evaluated for different environmental temperatures defining the consumption behavior for the large scale network. These scenarios cover the relevant mean daily temperatures $\{-3,3,7.5 \}\,\cels$. The interval spanned by $[-3,7.5]\, \cels$ exhibits a large optimization potential in terms of distributing the feed-in power. For colder or warmer environmental temperatures, either all or none of the energy capacities within the power plant will be used. The lower and upper temperature constraints are in line with standard operation conditions of heating networks. Table \ref{tab:testcase} presents a detailed description of the optimization scenarios under investigation.
\begin{table}[h!]
\centering
\begin{tabular}{@{}llllll@{}}
\hline
$\Tminhs,\Tmaxnet / \cels$ & $\minpreshous$, $\maxpreshous$, $\maxpressdiff/$bar  & $\pinlim$  & $\daytmp/\cels$ & $t_0,t_e/h$ & \dt/s \\ \hline
75, 110 & 3.5, 9.1, 2.5 & $ 0.5 (\consmax + \consmn)$ & $ -3, 3, 7.5$ & 0,72 & 300 \\
\hline
\end{tabular}
\caption{Description of the considered optimization scenarios. Test cases TC1-TC3 differ by the considered daily mean temperature changing from $-3\cels$ to $7.5\cels$. $\Tminhs,\Tmaxnet$ denote the temperature equivalents of the optimization constraints (\ref{eq_op_minTemdiff},\ref{eq_op_uptemp}).}
\label{tab:testcase}
\end{table}
The observation interval, in which the constraints and the objective function are evaluated, is set to $\dt = 300 \,$s, which is smaller than the typical plant operation interval of $900\,$s. This allows to approximate the underlying dynamics more precisely, while matching the relevant decision interval. The power constraint $\pinlim \in [\consmn, \consmax]$ is chosen within the mean ($\consmn$) and maximum ($\consmax$) daily consumption. While the mean consumption $\consmn$ naturally poses a lower limit for the maximum injected power, the maximum consumption $\consmax$ is an upper limit, since it can always be achieved by a stationary control. Due to the initialization with a stationary solution, during the first period the power restriction is relaxed to the maximum consumption. During this time, the output energy density and the corresponding constraints are shaped by the initial solution and not the control. The friction factor $\lambda$ entering the Darcy-Weisbach law (\ref{eq_consmom_sim}) is modeled by the Colebrook-White equation depending on the Reynolds number of the fluid, and the roughness and diameter of each pipeline. A computationally efficient approximation of the frictional model is achieved by setting a time-independent friction factor $\lambda_i, \; i \in \pipes$ for every pipeline by a-priori defining a suitable Reynolds number. The difference of maximal and minimal pressure levels at all consumption points in the flow network is given by $\maxpressdiff=$2.5\,bar and replaces the optimization constraints (\ref{eq_op_minPres_hs},\ref{eq_op_maxPres_hs}). For the gravitational constant we use a value of $g = 9.81\,$m\,s$^{-2}$. Specific heat capacity and density are set to $c_p = 4.16\,$kJ\,K$^{-1}$\,kg$^{-1}$, and $\rho = 1000\,$kg\,m$^{-3}$. To transform energy densities to temperatures presented in the section, we use the approximation 
\begin{align}
\edens \approx \rho c_p T.
\end{align}
This allows to transform the constraints for the minimum energy density at consumer stations (\ref{eq_op_minTemdiff}), as well as the maximum energy density in the network (\ref{eq_op_uptemp}) to their corresponding temperature values $\Tminhs$, and $\Tmaxnet$ displayed in tab.~\ref{tab:testcase}.\\

The power extraction $\cons$ is modeled using demand profiles typically employed in the simulation of heating networks \cite{SLP_BGW}. Specifically, each heat exchanger exhibits the power demand
\begin{align}
\cons_i(t) = c_i \cdot s_{m(i)}(t,\daytmp), \quad i \in \mathcal{H}.
\end{align}
The customer specific scaling factor $c_i$ represents an estimate for the total daily energy demand. The profile $s_m(t,\daytmp)$ models the time dependent demand for a given consumption class $m$ and the daily mean environmental temperature $\daytmp$. Every member of the class thus shares the same normalized profile while exhibiting an own specific consumption. The available demand profiles adjust the daily consumption by hourly scaling factors. The latter also depend on $\daytmp$, adjusting the relative weights of each hour in the daily consumption. Based on these hourly values, a spline interpolating the consumption is generated. In the considered network, the largest part of consumers belongs to the same consumption class $s_0$. To this end, all consumers are modeled by this class.\\

Since the typical consumption follows a periodic profile if the daily mean environmental temperature does not change, the control $\parcont$ is parameterized by a Fourier series,
\begin{align}
\parcont(t) = c_0 + \sum_{k=1}^K c_k \cos(k \omega t)+ \sum_{k=1}^K s_k \sin(k \omega t),
\end{align}
which approximates any control $u \in L_2$. Here, the frequency $\omega = 2 \pi/\timperiod$ is fixed to the period length $\timperiod$ of the consumption signal corresponding to 24\;h. The remaining Fourier coefficients act as parameters to be optimized.
\subsection{Optimal control for TC1}
The discussion of the numerical optimization results starts with TC1 simulating a mean daily temperature of $-3\cels$ with a mean and maximum consumption of $1.64\,$MW, and $2.29\,$MW.\\

\begin{figure}
\centering
\includegraphics[scale=1]{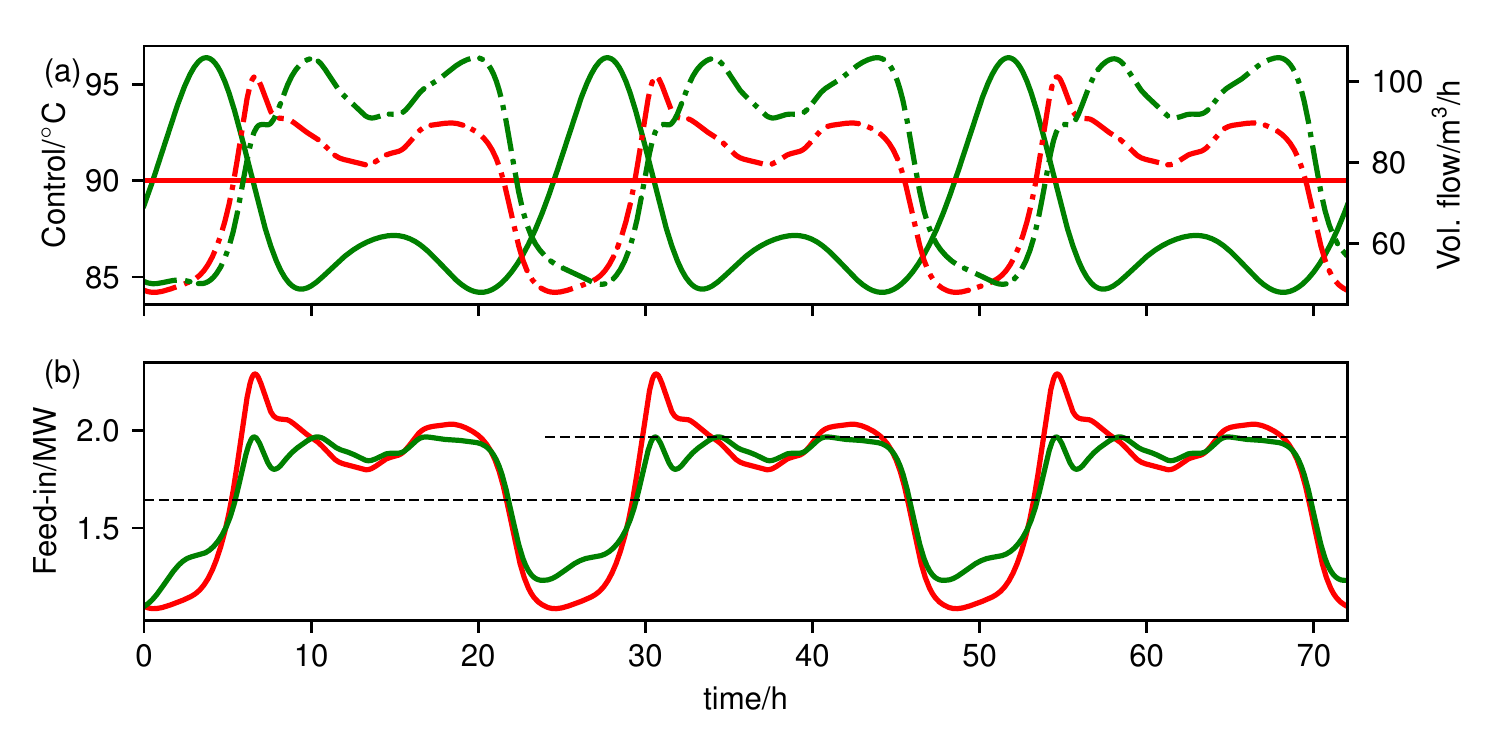}
\caption{Optimal control problem for TC1 at $-3\cels$ comparing the initial control (red) to the optimized control (green), obtained by the reduced model ROM1. Part (a) shows controls(solid lines) and the total volume flow injected at the power plant(dashed lines). Part (b) presents the feed-in power for both controls together with the mean consumption (lower dashed line) and the feed-in constraint $\pinlim$ (upper dashed line).}
\label{fig_conv_-3}
\end{figure}

\begin{figure}
\centering
\includegraphics[scale=1]{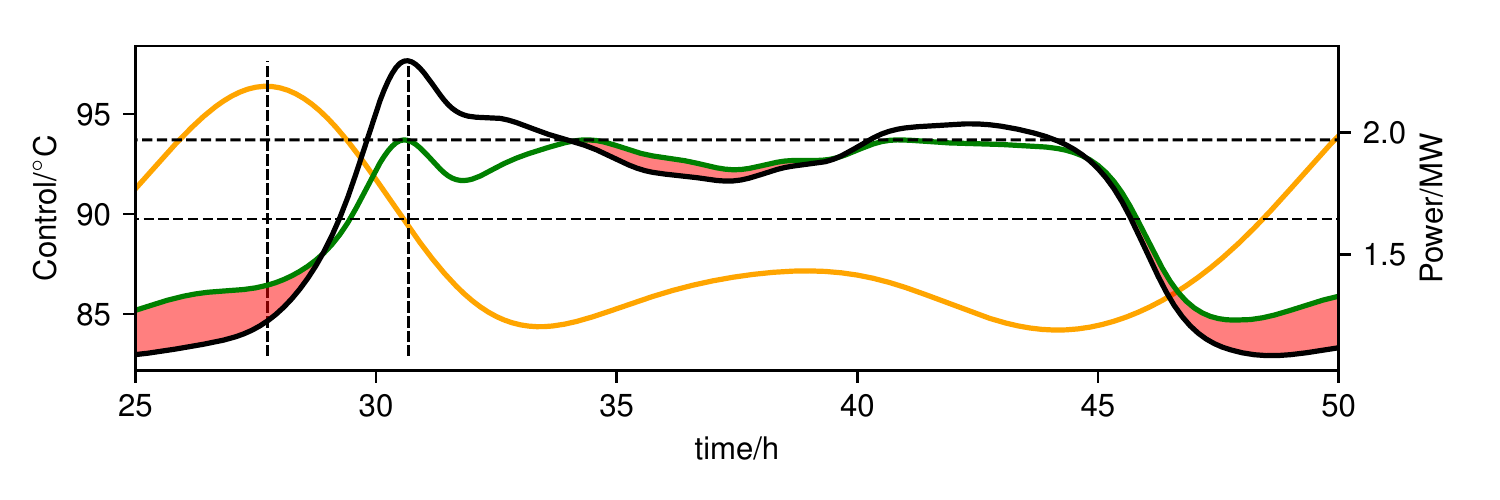}
\caption{Result of the optimal control problem for TC1 at $-3\cels$ visualizing the control at the power plant (orange) leading to a feed-in power (green) below the maximum constraint (upper dashed line). Red areas indicate regions in which pre-heating happens: The feed-in power exceeds the current consumption (solid, black line). Vertical, dashed lines visualize the time difference between the maximum injected temperature and the maximum consumption. The lower, dashed line indicates the mean consumption during one day as a guide for the eye.}
\label{fig_delay_-3}
\end{figure}

\begin{figure}
\centering
\includegraphics[scale=1]{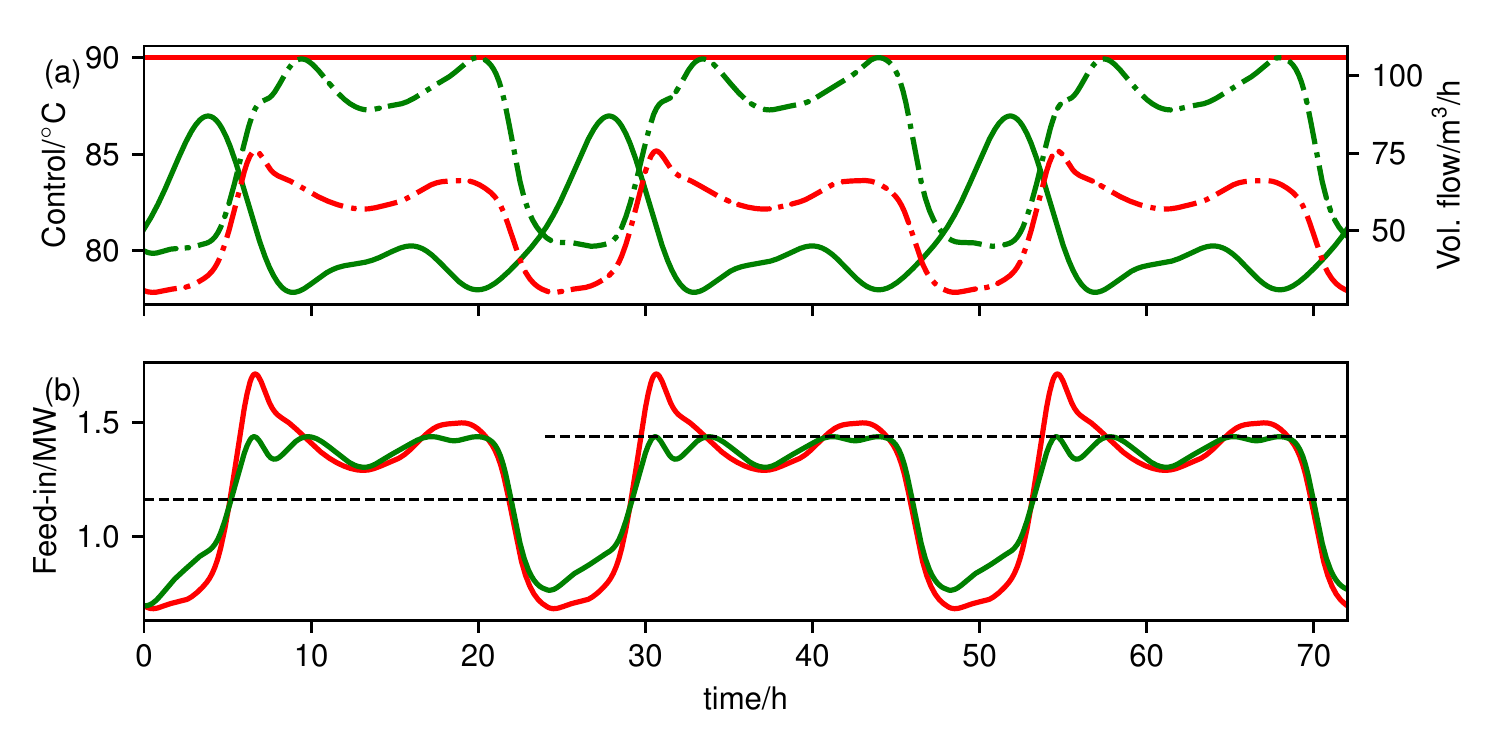}
\caption{Optimal control problem for TC2 at $3\cels$ comparing the initial control (red) to the optimized control (green), obtained by the reduced model ROM1. For a detailed explanation we refer to fig.~\ref{fig_conv_-3}.}
\label{fig_conv_3}
\end{figure}

Before comparing runtime and optimal controls obtained by different spatial discretizations, we analyze the optimal control suggested by ROM1. This model is gained by Galerkin projection of the upwind discretization FOM1, cf. tab.~\ref{tab:runtime_fullredTC1}. Based on a constant initial control of $90 \cels$, avoiding high feed-in peaks forces the control to increase its temporal variation, cf. fig.~\ref{fig_conv_-3}(a). The resulting feed-in power remains below the power constraint for all times cf. fig.~\ref{fig_conv_-3}(b). Since the mean temperature changes only slightly, the resulting total volume flow injected at the power plant also remains on the same level compared to the initial control, cf. fig.~\ref{fig_conv_-3}(a).
The suggested temporal variation of the thermal control which is necessary to limit the feed-in power induces pre-heating effects, visualized in fig.~\ref{fig_delay_-3}. Thus, reflecting the transport delay from source to the sink, in certain time intervals marked in red, the injected feed-in power exceeds the current consumption. More specifically, the maximum temperature level is injected at the power plant before the maximum aggregated consumption occurs at the consumption points. A common feature of the standardized consumption profiles are peak consumption around 6 a.m. and 6 p.m. Indeed, two pre-heating phases can be observed reflecting these two phases of high consumption, cf. fig.~\ref{fig_delay_-3}.

To avoid finite horizon effects, in which the suggested optimal control exploits the energy incorporated in the initial state, the setup is simulated for three periods of 24 hours. Focusing on the injected power coupling to volume flows at consumer stations as state variables, a state close to periodicity is reached quickly.
\begin{table}[h!]
\centering
\begin{tabular}{@{}l|llll@{}}
\hline
 & FOM 0 & FOM 1 & ROM 1 & \disref \\ \hline
runtime opt./s & 2765.0 & 2640.0 & 608.5 & 16100.0 \\ 
runtime sim./s & 2685.0 & 2585.0 & 557.0 & 15900.0 \\ 
\nDAEsolv & 16 & 9 & 11 & 12 \\ 
DOF & 775 & 1789 & 180 & 9538 \\ 
$J(u)$ & 55.9 &53.0 & 52.9 & 52.6 \\ 
$\violcon$ & 8.10$\times 10^{-3}$ & 9.48$\times 10^{-4}$ & 7.25$\times 10^{-4}$ & 0 \\ 
$\violfeed$ &5.41$\times 10^{-3}$ & 4.60$\times 10^{-4}$ & 4.49$\times 10^{-4}$ & 0 \\ 
$\maxreldev$ &3.03$\times 10^{-2}$ & 1.23$\times 10^{-3}$ & 5.28$\times 10^{-3}$ & 0 \\ 
\hline
\end{tabular}
\caption{Optimal controls and runtime comparison of TC1 ($-3 \,\cels$) for varying number of state variables (DOF) including full and reduced order models. The reference control $\tilde u$ results from the optimal control strategy described in section \ref{sec_det_optcont} using the reference discretization \disref. The feed-in $\feedin$ is measured by \disref\; based on the control suggested by each coarse model. The last row measures the relative error of outputs $y$ comparing reference model \disref, and each coarse model. }
\label{tab:runtime_fullredTC1}
\end{table}
\subsubsection*{Optimal controls resulting from different spatial discretizations}
The control resulting from ROM1 which is presented in fig.~\ref{fig_conv_-3} attains the feed-in constraint at several points of time. To check for feasibility of the solution, the suggested optimal controls resulting from both full- and reduced order models are compared and validated using a reference discretization \disref. The latter is given by an upwind discretization in space with a high number of finite volume cells. This will answer the question, whether a coarse, unreduced model is appropriate for an optimization task as well. To this end, table~\ref{tab:runtime_fullredTC1} compares the runtimes and approximation qualities of both full and reduced order models. While FOM0 denotes the minimal upwind discretization, in which each pipeline receives one finite volume cell, FOM1 is a finer discretization, with small approximation errors compared to \disref. Finally, ROM1 is the reduced order model obtained by reducing the discretized model FOM1. As quality indicators of the optimal control, the objective function, the relative $l_2$ error of the control and the maximal relative error of the feed-in constraint are considered. In addition, the maximum relative $l_2$ error of all outputs $y$ is measured. Here, $y^h_R \in \mathcal{H}$ refers to output $h$ determined by the reference discretization. Specifically, $\tilde u$ is the optimal control determined by the fine model \disref , and $\violfeed$ results from the control suggested by a coarse model simulated using \disref.\\

Focusing on TC1 and the suggested objective functions $J(u)$, FOM1, FOMR, and ROM1 converge to a comparable value with \disref\;taking the minimum of 52.6. FOM0 deviates clearly to 55.9. Regarding the relative deviation to the reference control, ROM1 shows the best approximation with $7.25 \times 10^{-4}$, followed by the full order models FOM1 and FOM0. The relative violation of the feed-in constraint is smallest for ROM1, followed by FOM1 and FOM0, while all models exhibit relative errors below one percent. Regarding the approximation quality of the outputs, FOM0 shows the expected strong diffusion, leading to a maximum relative error of $3.03 \times 10^{-2}$. Although this error does not affect the feasibility of the feed-in constraint, it leads to violations of the temperature constraints, measurable in practical applications. In contrast, ROM1 still approximates the outputs with an error of $5.28 \times 10^{-3}$.\\

Focusing on the runtimes for determination of the optimal control, ROM1 allows for a speed-up of 4.5 of the entire optimization compared to the coarsest and thus fastest possible unreduced model FOM0. In addition, the speed up compared to FOM1 amounts to 4.3, while achieving comparable results. The runtime of FOM0 results from a higher number of iterations necessary to fulfill the optimization tolerances. For the determination of the optimal control, resimulating the dynamics for a new control candidate takes the largest computational cost.

\subsubsection*{Different environmental temperatures}
Hereafter, test cases TC2, TC3 are discussed, simulating higher mean daily temperatures. In these scenarios, both mean and maximum power consumption decrease compared to TC1. To achieve comparable power constraints for different test cases, a relative constraint $\pinlim = \consmn + 0.5 (\consmax - \consmn)$ for the feed-in is chosen.\\

In contrast to TC1, the optimal control suggested for TC2 by ROM1 decreases the thermal control $\parcont$ compared to the initial control at $90\cels$. The corresponding average volume flow increases. This happens at the expense of increasing pumping costs which can be neglected as described in section \ref{sec_cont_prob}. Since the total power consumption for TC2 is smaller, the resulting increase in the volume flow does not exceed the level of TC1, allowing for identical pressure differences within the required constraints. In addition to a decreased mean value of the thermal control, also the temporal variation can be decreases. This leads to a smaller objective function of 14.6. Concerning error measures, the same observations discussed for TC1 apply. The speed-up of ROM1 compared to FOM0(FOM1) results in a factor of 4.2(5.4).\\

\begin{figure}[h!]
    \centering
    \begin{minipage}{.5\textwidth}
        \centering
        \includegraphics{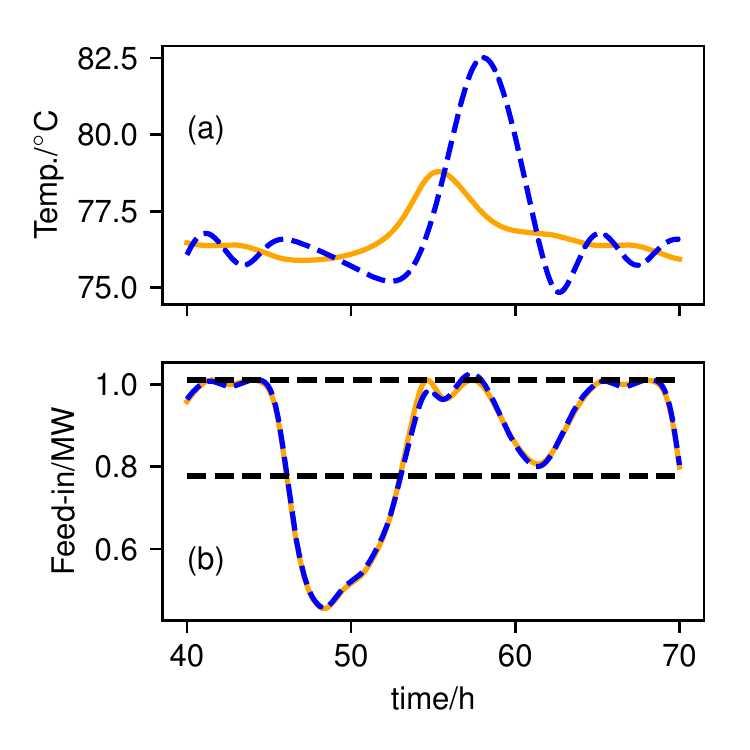}
    \end{minipage}%
    \begin{minipage}{0.5\textwidth}
        \centering
        \includegraphics{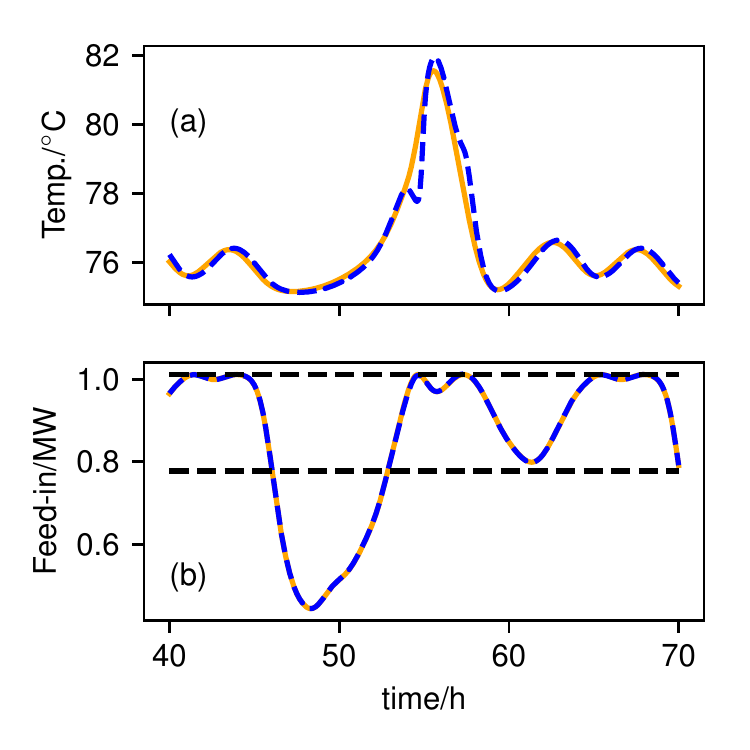}
    \end{minipage}
    \caption{Temperature signal for TC3 at the consumer exhibiting the largest relative $l_2$ error (top) and feed-in power (bottom) comparing FOM0 (left) and the reduced model ROM1 (right). The output of both models (orange, solid) is compared to their validation using the reference discretization \disref \; (blue, dashed line).}
        \label{fig:comp_state_approx}
\end{figure}

For the last scenario TC3 simulating a mean daily temperature of $7.5\cels$, the thermal control decreases in average value and temporal variation even further. The reduced consumption allows to increase the injected volume without violating pressure constraints. The injected flow temperature does now approaches the lower limit defined at $75\cels$. As observed for the other test cases, the approximation of the feed-in power is remarkably precise even for large deviations in the approximations of the outputs. For the coarsest discretization FOM0, the relative error of the feed-in constraint results in $1.46 \times 10^{-2}$. To illustrate this observation for TC3 , fig.~\ref{fig:comp_state_approx}(a) shows the output with the largest relative error compared to the reference discretization \disref. The robustness of the feed-in constraint towards errors in the state-space approximation is unexpected by (\ref{eq_feedin}), in which the feed-in depends on the volume flows at households defined by the thermal outputs $y$. Two explanations can be supplied for this effect. First, the feed-in depends on the sum of volume flows over all consumer points, allowing approximation errors to cancel. Furthermore, as the upwind discretization is a conservative finite volume scheme, the total stored energy is preserved on every discretization level.\\

In contrast to the large deviations observed for FOM0, ROM1 displayed in fig.~\ref{fig:comp_state_approx}(b) exhibits smaller errors in the output approximation. Specifically, it mainly deviates around the discontinuity at $t=54$\,h resulting from changes in the flux direction occurring in the dynamical simulation visible for the reference discretization. Again, the error in the feed-in constraint resulting from an imprecise approximation of flux changes is small by two reasons. First, only few consumers are affected by changing flux directions. Second, a change of flux directions is associated with the volume flow tending to zero, forcing the implied power to be zero as well. Hence, the absolute error in the power approximation remains small.
\begin{figure}
\centering
\includegraphics[scale=1]{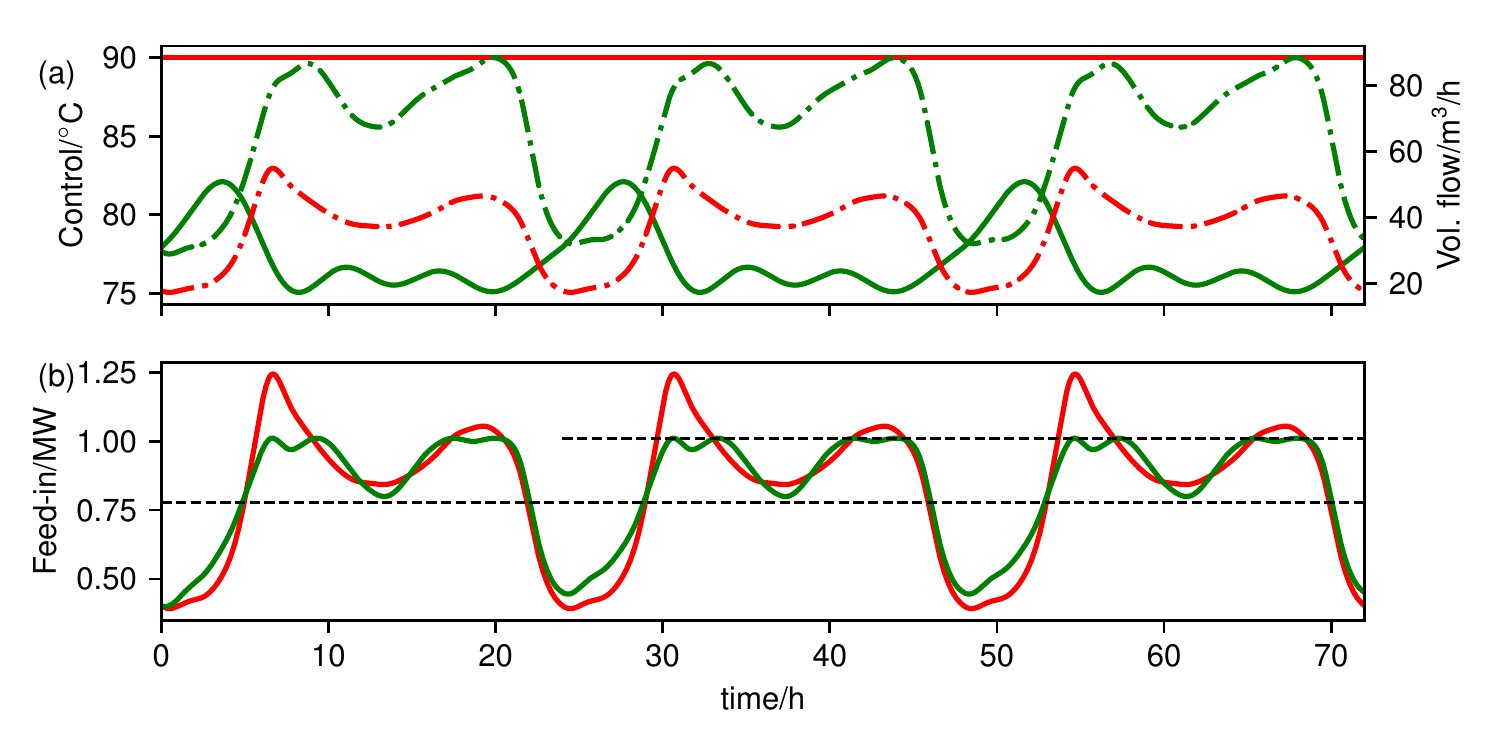}
\caption{Optimal control problem for TC3 at $7.5\cels$ comparing the initial control (red) to the optimized control (green), obtained by the reduced model ROM1. For a detailed explanation we refer to fig.~\ref{fig_conv_-3}. }
\label{fig_conv_7,5}
\end{figure}
\begin{table}[h!]
\centering
\begin{tabular}{@{}l|llll@{}}
\hline
 & FOM 0 & FOM 1 & ROM 1 & \disref \\ \hline
runtime opt./s & 2525.0 & 3210.0 & 596.2 & 18300.0 \\ 
runtime sim./s & 2455.0 & 3146.7 & 546.2 & 18100.0 \\ 
\nDAEsolv & 15 & 11 & 11 & 13 \\ 
DOF & 775 & 1789 & 180 & 9538 \\ 
$J(u)$ & 16.5 &14.8 & 14.8 & 14.6 \\ 
$\violcon$ & 7.10$\times 10^{-3}$ & 6.17$\times 10^{-4}$ & 5.26$\times 10^{-4}$ & 0 \\ 
$\violfeed$ &8.11$\times 10^{-3}$ & 6.86$\times 10^{-4}$ & 6.80$\times 10^{-4}$ & 0 \\ 
$\maxreldev$ &2.60$\times 10^{-2}$ & 1.25$\times 10^{-3}$ & 4.72$\times 10^{-3}$ & 0 \\ 
\hline
\end{tabular}
\caption{Optimization results and runtime comparison of TC 2($+3 \,\cels$) for different spatial discretizations including full and reduced order models. For a detailed explanation we refer the reader to table~\ref{tab:runtime_fullredTC1}.}
\label{tab:runtime_fullredTC2}
\end{table}
\begin{table}[h!]
\centering
\begin{tabular}{@{}l|llll@{}}
\hline
 & FOM 0 & FOM 1 & ROM 1 & \disref \\ \hline
runtime opt./s & 3165.0 & 5930.0 & 801.8 & 16400.0 \\ 
runtime sim./s & 3070.0 & 5795.0 & 735.2 & 16200.0 \\ 
\nDAEsolv & 19 & 21 & 15 & 13 \\ 
DOF & 775 & 1789 & 180 & 9538 \\ 
$J(u)$ & 6.7 &5.6 & 5.6 & 5.5 \\ 
$\violcon$ & 6.16$\times 10^{-3}$ & 6.14$\times 10^{-4}$ & 9.63$\times 10^{-4}$ & 0 \\ 
$\violfeed$ &1.46$\times 10^{-2}$ & 1.34$\times 10^{-3}$ & 1.37$\times 10^{-3}$ & 0 \\ 
$\maxreldev$ &2.35$\times 10^{-2}$ & 1.57$\times 10^{-3}$ & 5.01$\times 10^{-3}$ & 0 \\ 
\hline
\end{tabular}
\caption{Optimization results and runtime comparison of TC 3($+7.5 \,\cels$) for different spatial discretizations including full and reduced order models. For a detailed explanation we refer the reader to table~\ref{tab:runtime_fullredTC1}.}
\label{tab:runtime_fullredTC3}
\end{table}
\subsection{Speed-up of the reduced order model}
As discussed above, resimulating the dynamics is the central computational cost in determining an optimal control. One cause for the speed-up of the ROM is the few number of simulations of the forward problem to satisfy the tolerances for constraints and the objective function. The second cause is the speed-up resulting for a single solution of the DAE (\ref{eq_opt_DAE}), which is discussed subsequently. For the implicit midpoint rule with $n_t$ time steps, the computational cost $\costsim$ splits into the following parts,

\begin{align}
\costsim = n_t \cdot (n_{imp} \cdot(c_h + c_f)  + c_J),
\end{align}

where $n_{imp}$ is the number of iterations to solve for the upcoming timestep, $c_h$ is the cost of solving the algebraic equations (\ref{eq_cons_vol}-\ref{eq_hous}) defining the flow field, and $c_f$ is the cost of evaluating the differential part of the DAE. Finally, $c_J$ denotes the cost to determine the Jacobian of the DAE with respect to energy densities. Since only the thermal transport is reduced, $c_h$ is identical for both full and reduced order models. The cost for the evaluation of the ODE scales with the number of entries implied in the system operators $A(v), A_r(v)$ defined in (\ref{eq_sys},\ref{eq_sys_red}). Using the upwind discretization, $\elrmat{n}{n}{A(v)}$ is a sparse matrix. For typical networks, the number of nonzero entries can be limited by $3n$. In contrast, the system operator $\elrmat{r}{r}{A_r(v)}$ resulting from a Galerkin projection is dense with $r^2$ nonzero entries. As a consequence, the number of reduced states needs to be significantly smaller to reduce the computational cost. For the coarse discretization sufficient for the determination of an optimal control, this degree of reduction is barely possible. The key saving in applying the reduced order model stems from the computation of the Jacobian matrix. Based on the system operator description and the definition of the DAE (\ref{eq_sys})
\begin{align*}
A(v) &= \sum_{i=1}^{\nvflow} \wgtfun^v_i(v)  A^v_i,
\end{align*}
the Jacobian reads
\begin{align}
\jacob_{\dot e}(v,e) = \frac{\partial \dot e}{\partial e} = A(v)+ \sum_{i=1}^{\nvflow} \frac{\partial \wgtfun^v_i(v)}{\partial e}  (A^v_i e + B^v_i \uT).
\end{align}
Since by the algebraic equation (\ref{eq_hous}) each velocity depends on the energy densities at households, $\jacob_{\dot e}(v,e)$ carries significantly more non-zero entries in the full order case than $A(v)$, cf. tab.~\ref{tab:populaton_jac}. In contrast, since the reduced order operator $A_r(v)$ is already densely populated by Galerkin projection, the number of non-zero entries does not increase for the reduced Jacobian.
\begin{table}[h!]
\centering
\begin{tabular}{@{}l|llll@{}}
\hline
 & FOM 0 & FOM 1 & ROM 1 & \disref \\ \hline
\# nonzero entries $A(v)$ & 1555 & 3583 & 6827 & 19081 \\ 
\# nonzero entries $\jacob_{\dot e}(v,e)$ & 67837 &139922 & 11211 & 662200\\
DOF & 775 & 1789 & 180 & 9538 \\ 
\hline
\end{tabular}
\caption{Maximal population density determined in the simulation of system operator and Jacobian for different discretizations. DOF denotes the number of differential state variables.}
\label{tab:populaton_jac}
\end{table}
\section{Conclusions}
In this contribution, we discussed the optimal control of district heating networks utilizing a reduced order model (ROM). The suggested optimal controls resulting from minimizing the temporal variation of the control successfully limit the maximum feed-in power to the average of mean and maximum total consumption. In addition, practically relevant constraints on temperature and pressure are included reproducing realistic operation conditions. For the presented scenarios, this allows to avoid the usage of additional, unfavorable sources of energy. The quadratic in-state DAE of index 1 is split into a thermal and a flow defining part, allowing to describe the thermal part as a linear, parameter-varying system with the velocity field acting as the parameter. Using a greedy strategy, relevant velocity configurations are implied in a global Galerkin projection forming a Lyapunov stable ROM. While the ROM approximates both relevant state variables and gradient information sufficiently fine for the determination of an optimal control, it speeds up the entire optimization phase by at minimum a factor 4, compared to even coarse levels of upwind discretizations used as full order models. For distinct test scenarios we observe even higher speed-ups of 7.3. Thus, the ROM gaps the bridge towards the determination of an optimal control within an online planning. The effectiveness of the ROM is demonstrated for an existing large scale network in which different pipelines change their flux direction dynamically. Runtime and approximation quality are studied for multiple real world scenarios including varying daily mean temperatures. This allows to apply the presented model to other networks and operation conditions relevant in practice.\\

For further research, we study the benefits of the reduced order model in a feedback control, in which additional advantages might result from its reduced state space dimension. Furthermore, a comparison to direct optimization approaches will be interesting, in which the transport dynamics of the network directly appear as optimization constraints, avoiding to resimulate the network dynamics in each iteration.

\section*{Acknowledgments}
We acknowledge the financial support by the German Federal Ministry BMBF within the project \textit{EiFer} (F\"orderkennzeichen: 05M18AMB), and by the German Federal Ministry BMWi within the project \textit{DYNEEF} (F\"orderkennzeichen: 03ET1346B).

\end{document}